\documentclass[11pt]{amsart2000}
\usepackage{amssymb}
\newtheorem{theorem}{Theorem}[section]
\newtheorem{corollary}[theorem]{Corollary}
\newtheorem{lemma}[theorem]{Lemma}
\newtheorem{proposition}[theorem]{Propostion}
\theoremstyle{definition}
\newtheorem{example}[theorem]{Example}
\newtheorem{question}[theorem]{Question}
\newtheorem{definition}[theorem]{Definition}

\input xypic
\def\cl{\operatorname{cl}}
\def\cli#1{\operatorname{cl}_{#1}}

\begin{document}
\setlength{\unitlength}{0.01in}
\linethickness{0.01in}
\begin{center}
\begin{picture}(474,66)(0,0)
\multiput(0,66)(1,0){40}{\line(0,-1){24}}
\multiput(43,65)(1,-1){24}{\line(0,-1){40}}
\multiput(1,39)(1,-1){40}{\line(1,0){24}}
\multiput(70,2)(1,1){24}{\line(0,1){40}}
\multiput(72,0)(1,1){24}{\line(1,0){40}}
\multiput(97,66)(1,0){40}{\line(0,-1){40}}
\put(143,66){\makebox(0,0)[tl]{\footnotesize Proceedings of the Ninth Prague Topological Symposium}}
\put(143,50){\makebox(0,0)[tl]{\footnotesize Contributed papers from the symposium held in}}
\put(143,34){\makebox(0,0)[tl]{\footnotesize Prague, Czech Republic, August 19--25, 2001}}
\end{picture}
\end{center}
\vspace{0.25in}
\setcounter{page}{171}
\title{Sequence of dualizations of topological spaces is finite}
\author{Martin Maria Kov\'ar}
\address{Department of Mathematics\\
Faculty of Electrical Engineering and Computer Science\\
Technical University of Brno\\
Technick\'a 8\\
616 69 Brno, Czech Republic}
\email{kovar@dmat.fee.vutbr.cz}
\begin{abstract}
Problem 540 of J. D. Lawson and M. Mislove in Open Problems in Topology
asks whether the process of taking duals terminate after finitely many
steps with topologies that are duals of each other. 
The problem for $T_1$ spaces was already solved by G. E. Strecker in 1966. 
For certain topologies on hyperspaces (which are not necessarily $T_1$),
the main question was in the positive answered by Bruce S. Burdick and
his solution was presented on The First Turkish International Conference
on Topology in Istanbul in 2000. 
In this paper we bring a complete and positive solution of the problem for
all topological spaces. 
We show that for any topological space $(X,\tau)$ it follows
$\tau^{dd}=\tau^{dddd}$.
Further, we classify topological spaces with respect to the number of
generated topologies by the process of taking duals.
\end{abstract}
\thanks{Reprinted from 
Theoretical Computer Science B,
in press,
Martin Maria Kov\'ar,
The solution to {P}roblem 540,
Copyright (2002),
with permission from Elsevier Science \cite{MMK}.}
\thanks{This research is supported by the grant no. 201/00/1466 of the
Grant Agency of the Czech Republic and by the research intention of the
Ministry of Education of the Czech Republic CEZ: J22/98:262200012}
\thanks{Martin Maria Kov\'ar,
{\em Sequence of dualizations of topological spaces is finite},
Proceedings of the Ninth Prague Topological Symposium, (Prague, 2001),
pp.~171--179, Topology Atlas, Toronto, 2002}
\maketitle

\setcounter{section}{-1}
\section{Definitions and Notation}

Through the paper we mostly use the standard topological notions. 
The main source of definitions of some newer concepts is the paper 
\cite{LM} or the book \cite{Vi} in which the reader can find many 
interesting connections to modern, computer science oriented topology. 
In this paper, all topological spaces are assumed without any separation
axiom in general.
Let $(X,\tau)$ be a topological space. 
Recall that the preorder of specialization is a reflexive and transitive 
binary relation on $X$ defined by 
$x\leqslant y$ if and only if $x\in \cl\left\{y\right\}$. 
This relation is antisymmetric and hence a partial ordering
if and only if $X$ is a $T_0$ space. 
For any set $A\subseteq X$ we denote 
$${\uparrow A}= \left\{x : x\geqslant y \text{ for some } y\in A
\right\}$$
and 
$${\downarrow A}= \left\{x : x\leqslant y \text{ for some } y\in A 
\right\}.$$
In particular, for a singleton $x\in X$ it is clear that 
${\downarrow\left\{x\right\}} = \cl\left\{x\right\}$. 
A set is said to be saturated in $(X,\tau)$ if it is the intersection of 
open sets. 
One can easily check that a set $A$ is saturated in $(X,\tau)$ if and 
only if $A={\uparrow A}$, that is, if and only if the set $A$ is closed 
upward with respect to the preorder of specialization of $(X, \tau)$. 
By the dual topology $\tau^d$ for a topological space $(X,\tau)$ we mean
the topology on $X$ generated by taking the compact saturated sets of
$(X,\tau)$ as a closed base. 
It is not very difficult to show that the dual operator switches the 
preorder of specialization, that is, $x\leqslant y$ in $(X,\tau)$ if 
and only if $x\geqslant y$ in $(X,\tau^d)$. 
Let $\psi$ be a family of sets. 
We say that $\psi$ has the finite intersection property, or briefly, that
$\psi$ has f.i.p., if 
$$\mbox{
for every $P_1, P_2,\dots, P_k\in\psi$ it follows
$P_1\cap P_2\cap \dots\cap P_k\ne\varnothing$.
}$$
For the definition of compactness, we do not assume any separation axiom.
We say that a subset $S$ of a topological space $(X,\tau)$ is compact if
every open cover of $S$ has a finite subcover. 
Let $(X,\tau)$ have a closed base $\Phi$. 
It is well-known that $S\subseteq X$ is compact if and only if for 
every family $\zeta\subseteq\Phi$ such that the family 
$\left\{S\right\}\cup\zeta$ has f.i.p.\ it follows 
$S\cap(\bigcap\zeta)\ne\varnothing$, or equivalently, if and only if every 
filter base $\eta\subseteq\Phi$ such that every element of $\eta$ meets 
$S$ the filter base $\eta$ has a cluster point in $S$.

\section{Introduction}

Problem 540 of J. D. Lawson and M. Mislove \cite{LM} in Open Problems in 
Topology (J. van Mill, G. M. Reed, eds., 1990) asks 
\begin{itemize}
\item 
which topologies can arise as dual topologies and
\item 
whether the process of taking duals terminate after finitely many steps
with the topologies that are duals of each other. 
\end{itemize}

For $T_1$ spaces, the solution of (2) follows simply from the fact that in
$T_1$ spaces every set is saturated and hence the dual operator $d$
coincide with the well-known compactness operator $\rho $ of J. de Groot,
G. E. Strecker and E. Wattel \cite{GSW}. 
For more general spaces, the question (2) was partially answered by Bruce
S. Burdick who found certain classes of (in general, non-$T_1$) spaces for
which the process of taking duals of a topological space $(X,\tau)$
terminates by $\tau^{dd}=\tau^{dddd}$ --- the lower Vietoris topology on
any hyperspace, the Scott topology for reverse inclusion on any
hyperspace, and the upper Vietoris topology on the hyperspace of a regular
space. 
B. Burdick presented his paper on The First Turkish International
Conference on Topology in Istanbul 2000 \cite{Bu}. 
In this paper we will show that an analogous result, that is,
$\tau^{dd}=\tau^{dddd}$ holds for arbitrary and general topological spaces
without any restriction.

\section{Main Results}

Through this paper, let $\Phi$ denote a certain closed topology base of a
topological space $(X,\tau)$. 
Let $\Phi^d$ be the collection of all compact saturated sets in $(X,\tau)$
and hence the closed topology base of $(X,\tau^d)$. 
For $n=2, 3 \dots$ let $\Phi^{d^n}$ be the collection of all compact
saturated sets in $(X,\tau^{d^{n-1}})$ and hence the closed topology base
of $(X,\tau^{d^n})$. 
The following lemma is due to Bruce S. Burdick \cite{Bu}. 
We repeat it because of completeness.

\begin{lemma}[B. S. Burdick, 2000)]
Let $C$ be compact in $(X,\tau)$ and $P \in\Phi^{dd}$. 
Then $C\cap P$ is compact in $(X,\tau)$.
\end{lemma}

\begin{proof} 
Let $\zeta\subseteq\Phi$ be a filter base such that its every element 
meets $C\cap P$. 
Let $$\eta=\left\{ {\uparrow(C\cap D)} : D\in\zeta\right\}.$$
Then $\eta\subseteq\Phi^d$ and $\left\{P\right\}\cup\eta$ has f.i.p. 
Since $P \in\Phi^{dd}$ it follows that there exists some 
$$x\in P\cap(\bigcap\eta)\ne \varnothing.$$
Hence, $x\in P\cap {\uparrow(C\cap D)}$ for every $D\in\zeta$,
which implies that there is some $t\in C\cap D$ (depending on the choice
of $D$) such that $x\geqslant t$. 
Then 
$${\downarrow \left\{x\right\}} \cap C\cap D\ne \varnothing
\quad
\mbox{for every $D\in\zeta$},$$ 
so the collection 
$$
\left\{C\right\}\cup (\left\{{\downarrow\left\{x\right\}}\right\}\cup\zeta)
$$ 
has f.i.p.
Since $C$ is compact and $\downarrow \left\{x\right\}$ as well as every 
element of $\zeta$ is closed, there is some 
$$
y \in 
C \cap {\downarrow \left\{x \right\}} \cap (\bigcap\zeta)
\subseteq C \cap P \cap (\bigcap\zeta).
$$ 
Hence, $y$ is a cluster point of the closed filter base $\zeta$ in
$C\cap P$.
\end{proof}

\begin{lemma} 
Let $C\in\Phi^d$, $M\in\psi\subseteq\Phi^{dd}$ and let 
$\left\{C\right\}\cup\psi$ has f.i.p. 
Then there exist $\xi(M)\in M$ such that
$$
\left\{C\right\}\cup\psi\cup\left\{{\downarrow\left\{\xi(M)\right\}}\right\}
$$ 
has f.i.p.
\end{lemma}

\begin{proof} 
Let 
$$
\varphi = 
\left\{ {\uparrow (C\cap P_1\cap P_2\cap \dots\cap P_k)} : 
P_1, P_2,\dots,P_k\in\psi\right\}.
$$
From Lemma~2.1 it follows that $\varphi\subseteq\Phi^d$ and 
$\left\{M\right\}\cup\varphi$ has f.i.p. 
Hence, since $M\in\Phi^{dd}$ is compact in $(X,\tau^d)$, there exist 
$$\xi(M)\in M \cap(\bigcap\varphi).$$
Then for every $P_1, P_2,\dots,P_k\in\psi$ it follows that 
$$\xi(M)\in {\uparrow (C\cap P_1\cap P_2\cap \dots\cap P_k)}$$ 
so there exist $t\in C\cap P_1\cap P_2\cap \dots\cap P_k$ (depending on
the choice of $P_1, P_2,\dots, P_k$) with $\xi(M)\geqslant t$. 
Then $t\in{\downarrow\left\{\xi(M)\right\}}$ which implies that
$$
{\downarrow\left\{\xi(M)\right\}} \cap C\cap P_1\cap P_2\cap \dots\cap P_k
\ne\varnothing.
$$
It follows that 
$$\left\{C\right\} \cup \psi \cup 
\left\{{\downarrow\left\{\xi(M)\right\}}\right\}$$ 
has f.i.p.
\end{proof}

\begin{lemma}
Let $C\in\Phi^d$, $\psi\subseteq\Phi^{dd}$ and let 
$\left\{C\right\}\cup\psi$ has f.i.p. 
Then for every $M\in\psi$ there exist $\xi(M)\in M$ such that
$$
\left\{C\right\}\cup\left\{{\downarrow\left\{\xi(M)\right\}} : 
M\in\psi\right\}
$$ 
has f.i.p.
\end{lemma}

\begin{proof}
Let $\psi=\left\{M_\alpha : \alpha<\mu\right\}$ where $\mu$ is an
ordinal number. 
From Lemma~2.2 it follows that there exist $\xi(M_0)\in M_0$ such that
$$
\left\{C\right\}\cup\psi\cup\left\{{\downarrow\left\{\xi(M_0)\right\}}\right\}
$$ 
has f.i.p. 
Suppose that for some $\beta<\mu$ and every $\alpha<\beta$ there exist
$\xi(M_\alpha)\in M_\alpha$ such that, in the notation
$$
\chi_\alpha = \psi \cup \left\{{\downarrow \left\{\xi(M_\gamma)\right\}} : 
\gamma\leq\alpha\right\},
$$ 
the family $\left\{C\right\}\cup \chi_\alpha$ has f.i.p. 
Let $\chi=\bigcup_{\alpha<\beta}\chi_\alpha$. 
Obviously, the family $\left\{C\right\}\cup \chi$ has f.i.p. and
$M_\beta\in\psi\subseteq\chi \subseteq\Phi^{dd}$. 
Then, by Lemma~2.2 there exist $\xi(M_\beta)\in M_\beta$ such that
$$
\left\{C\right\} \cup \chi \cup 
\left\{{\downarrow\left\{\xi(M_\beta)\right\}}\right\}
$$ 
has f.i.p. 
But 
$$\chi_\beta =\chi \cup 
\left\{{\downarrow\left\{\xi(M_\beta)\right\}}\right\}$$
which implies that the family $\left\{C\right\}\cup \chi_\beta$ has
f.i.p. 
By induction, we have defined $\xi(M_\beta)\in M_\beta$ for every
$\beta<\mu$. 
Obviously, the family
$$
\left\{C\right\} \cup (\bigcup_{\beta<\mu}\chi_\beta) =
\left\{C\right\} \cup \psi \cup
\left\{{\downarrow\left\{\xi(M_\beta)\right\}} : \beta<\mu\right\}
$$ 
has f.i.p which implies that also its subfamily 
$$
\left\{C\right\} \cup 
\left\{{\downarrow\left\{\xi(M)\right\}} : M\in\psi\right\}
$$ 
has f.i.p. 
\end{proof}

\begin{theorem}
For any topological space $(X, \tau)$ it follows 
$\Phi^d \subseteq\Phi^{ddd}$.
\end{theorem}

\begin{proof}
Let $C\in\Phi^d$ and let $\psi\subseteq\Phi^{dd}$ be a family such that
$\left\{C\right\}\cup\psi$ has f.i.p. 
From Lemma~2.3. it follows that for every $M\in\psi$ there exist
$\xi(M)\in M$ such that 
$$\left\{C\right\} \cup \left\{{\downarrow\left\{\xi(M)\right\}} : 
M\in\psi\right\}$$ 
has f.i.p.
But $C$ is compact in $(X,\tau)$ and for every $M\in\psi$, the sets 
$${\downarrow \left\{\xi(M)\right\}} = 
\cl\left\{\xi(M)\right\} \subseteq 
M$$
are closed. 
Hence 
$$
\varnothing \ne 
C \cap (\bigcap\left\{{\downarrow \left\{\xi(M)\right\}} : 
M\in\psi\right\}) \subseteq 
C\cap(\bigcap\psi).
$$
It follows that $C$ is compact in $(X,\tau^{dd})$. But $C$ is saturated in
$(X,\tau)$ which is the same as $C$ is saturated in $(X,\tau^{dd})$. 
Hence $C\in\Phi^{ddd}$. 
\end{proof}

\begin{corollary} 
For any topological space $(X, \tau)$ it follows $\Phi^{dd} =\Phi^{dddd}$.
\end{corollary}

\begin{proof}
From Theorem 2.4 it follows that $\Phi^d\subseteq\Phi^{ddd}$
which implies that $\Phi^{dd}\supseteq\Phi^{dddd}$. 
Applying Theorem 2.4 to the space $(X,\tau^d)$ we obtain that
$\Phi^{dd}\subseteq\Phi^{dddd}$.
\end{proof}

\begin{corollary}
For any topological space $(X, \tau)$ it follows $\tau^{dd} =\tau^{dddd}$.
\end{corollary}

\section{Some Classification Notes}

Now, let $e$ be an identity operator on the class of all topological
spaces, i.e. $e(\tau)=\tau$ for a topology $\tau$ of a space $(X,\tau)$.
Then $e$, $d$, $dd$ and $ddd$ form, with the composition operation
$\circ$, a commutative monoid having the unit element $e$ and the
following multiplicative table:
$$
\begin{array}{|r|l|l|l|l|}
\hline
\circ	&e	&d	&dd	&ddd	\\
\hline
e	&e	&d	&dd	&ddd	\\
\hline
d	&d	&dd	&ddd	&dd	\\
\hline
dd	&dd	&ddd	&dd	&ddd	\\
\hline
ddd	&ddd	&dd	&ddd	&dd	\\
\hline
\end{array}
$$
\centerline{\bf Table 1.}
\smallskip

Then it is natural to introduce the following classes 
$G_1$, 
$G_2 = G_{2_a}\cup G_{2_b}$, 
$G_3 = G_{3_a}\cup G_{3_b}\cup G_{3_c}$ and
$G_{4}$ of topological spaces, where: 
$$G_1 = \left\{(X,\tau) : \tau^d=\tau\right\},$$
$$
G_{2_a} = \left\{(X,\tau) : \tau^{dd}=\tau\right\},\ 
G_{2_b} = \left\{(X,\tau) : \tau^{dd}=\tau^d\right\},
$$
$$
\begin{array}{l}
G_{3_a} = \left\{(X,\tau) : \tau^{ddd}=\tau\right\},\\
G_{3_b} = \left\{(X,\tau) : \tau^{ddd}=\tau^d\right\},\\ 
G_{3_c} = \left\{(X,\tau) : \tau^{ddd}=\tau^{d}\right\},
\end{array}
$$
$$G_4 = \left\{(X,\tau) : \tau^{dddd} = \tau^{dd}\right\} = \mathsf{Top}.$$

Obviously, each class $G_n$ consists exactly of those topological spaces
such that the process of taking duals generates at most $n=1, 2, 3, 4$
different topologies while the class $G_4$ contains all topological
spaces. Hence, we suggest the following terminology:

\begin{definition}
We say that a topological space $(X,\tau)$ is {\it $n$-generative} if
$n=1,2,3,4$ is the least number such that $(X,\tau)$ belongs to the class
$G_n$.
\end{definition}

\begin{proposition}
The classes $G_n$ satisfy the following relationships:
\begin{itemize}
\item[(i)] $G_1 =G_{3_a}= G_{2_a}\cap G_{2_b}=G_{2_a}\cap G_{3_c}$
\item[(ii)] $G_{2_b}= G_{3_b}\cap G_{3_c}$
\item[(iii)] $G_3=G_{3_b}\cup G_{3_c}$
\end{itemize}
\end{proposition}

The proof is an easy consequence of the definitions of the classes $G_n$
and the identity $dd=dddd$ of Corollary 2.6. 
We leave it to the reader.
Hence, the relationships between the classes $G_n$ can be described by the
following diagram:
$$
\diagram
		&G_4\drline\dlline	&		\\
G_{3_b}\dline	&		 	&G_{3_c}\dline	\\
G_{2_a}		&			&G_{2_b}\ullline\\
		&G_1\urline\ulline	&
\enddiagram
$$
\smallskip
\centerline{\bf Diagram 1.}
\smallskip

Since for $T_1$ spaces the dual operator $d$ coincide with the compactness
operator $\rho$ of J. de Groot, G. E. Strecker and E. Wattel \cite{GSW} we
immediately obtain that 
$G_1\subsetneq G_{2_a}\subsetneq G_{3_b}\subsetneq G_4$, 
$G_1\subsetneq G_{2_b}\subsetneq G_{3_b}$ and $G_{3_c}\subsetneq G_4$ 
(the paper \cite{GHSW} is a good source of proper examples of such spaces,
see also the section 4). 
Hence it remains to determine whether the class $G_{2_b}$ is a proper
subclass of the class $G_{3_c}$.

\begin{lemma}
Let $\tau=\tau^d$ for a topological space $(X,\tau)$.
Then for every subset $A\subseteq X$ it follows that 
${\uparrow A} = {\downarrow A}$.
\end{lemma}

\begin{proof}
Let $x\in X$. 
Then 
$\downarrow\left\{x\right\} = \cli \tau \left\{x\right\}$ and
${\uparrow \left\{x\right\}} = \cli {\tau ^d} \left\{x\right\}$. 
Since $\tau=\tau^d$ it follows that 
${\uparrow \left\{x\right\}} = \downarrow\left\{x\right\}$. 
Then 
$$
\begin{array}{lll}
{\uparrow A}&
=&
\bigcup_{x\in A} {\uparrow \left\{x\right\}}\\
&
=& 
\bigcup_{x\in A}\downarrow\left\{x\right\}\\
& 
=&
\downarrow A.
\end{array}
$$
\end{proof}

\begin{lemma}
Let ${\uparrow A} = {\downarrow A}$ for every subset $A\subseteq X$. 
If $S\in\Phi^{dd}$ and $H\in\Phi$ then $S\cap H \in \Phi^{dd}$.
\end{lemma}

\begin{proof}
Take $\zeta\subseteq\Phi^d$ such that the family 
$\left\{S\cap H \right\}\cup\zeta$ has f.i.p. 
We put 
$$\eta=\left\{{\uparrow(H\cap D)} : D\in\zeta \right\}.$$
It follows that $\eta\in\Phi^d$ and the family $\left\{S\right\}\cup \eta$
has f.i.p. 
Since $S\in\Phi^{dd}$, we have 
$$S\cap(\bigcap\eta)\ne\varnothing.$$
It follows that there exist some $x\in S\cap(\bigcap\eta)$. 
Then $x\in {\uparrow (H\cap D)}$ for every $D\in\zeta$. 
Hence, there exist $t\in H\cap D$ (depending on the choose of $D$) such
that $x\geqslant t$. 
Obviously, $x\in D$ for every $D\in\zeta \subseteq\Phi^d$ since the
elements of $\Phi^d$ are saturated and hence closed upward with respect
to the specialization preorder. 
Since $H$ is a closed set, it follows that 
$$H = {\downarrow H} = {\uparrow H}$$
which implies that $x\in H$. Therefore, we have 
$$x\in S\cap H\cap(\bigcap\zeta)$$ 
and so $S\cap H \in \Phi^{dd}$.
\end{proof}

\begin{proposition}
The equality $G_{2_b}=G_{3_c}$ holds.
\end{proposition}

\begin{proof}
Let $(X,\tau)$ be a space that belongs to the class $G_{3_c}$. 
Then $\tau^{dd}=\tau^{ddd}$. 
Applying Lemma 3.3 to the space $(X,\tau^{dd})$ we obtain that for any
subset $A\subseteq X$ it follows that ${\uparrow A}={\downarrow A}$. 
We will show that $(X,\tau)$ belongs to the class $G_{3_b}$. 
It follows from Theorem 2.4 that $\Phi^d\subseteq\Phi^{ddd}$.
Conversely, let $S\in\Phi^{ddd}$. 
It means that $S$ is compact saturated in $(X,\tau^{dd})$ which is the
same as in $(X,\tau^{ddd})$. 
Hence $S\in\Phi^{dddd}=\Phi^{dd}$ by Corollary 2.5. 
We will show that $S\in\Phi^d$. 
Let $\zeta\subseteq\Phi$ be a family such that $\left\{S\right\}\cup\zeta$
has f.i.p. 
We put 
$$\eta=\left\{S\cap H : H\in \zeta\right\}.$$
It follows from Lemma 3.4 that $\eta\subseteq\Phi^{dd}$ and, obviously,
the family $\left\{S\right\}\cup \eta$ has f.i.p. 
Then 
$$
S\cap(\bigcap\zeta)=S\cap(\bigcap\eta)\ne\varnothing,
$$
which implies
that $S\in\Phi^d$. 
Then $\Phi^d=\Phi^{ddd}$ and so $\tau^d=\tau^{ddd}$. 
Hence, the space $(X,\tau)$ belongs to the class $G_{3_b}$. 
Now, we have that $G_{3_c}\subseteq G_{3_b}$ and 
$G_{2_b} = G_{3_b}\cap G_{3_c}$ by Proposition 3.2. 
Together it gives that $G_{2_b}=G_{3_c}$.
\end{proof}

\begin{corollary} 
The equality $G_3=G_{3_b}$ holds.
\end{corollary} 

Now, we can slightly improve Diagram 1:
$$
\diagram
	&G_4\dline		&	\\ 
	&G_3\dlline\drline	&	\\
G_{2_a}	&			&G_{2_b}\\
	&G_1\urline\ulline	&
\enddiagram
$$
\centerline{\bf Diagram 2.}
\smallskip

The two following simple results give a partial answer to the question
which topologies can arise by the process of taking duals.

\begin{corollary} 
Let $(X,\tau)$ be a topological space such that $\tau$ arise as a dual
topology. 
Then $(X,\tau)$ belongs to the class $G_3$ and does not belong to the
class $G_{2_b}\smallsetminus G_1$.
\end{corollary}

\begin{proof} 
If $\tau=\sigma^d$ where $\sigma$ is a topology on $X$, then 
$$\tau^d=\sigma^{dd}=\sigma^{dddd}=\tau^{ddd}.$$
Hence $(X,\tau)$ belongs to the class $G_3$. 
Suppose that $(X,\tau)$ belongs to the class $G_{2_b}$. 
Then 
$$\sigma^{dd}=\tau^d=\tau^{dd}=\sigma^{ddd}$$ 
which implies that 
$(X,\sigma)$ belongs to the class $G_{3_c}$. 
By Proposition 3.5 it follows that $G_{3_c}=G_{2_b}$ and so 
$$\tau=\sigma^d=\sigma^{dd}=\tau^d.$$
Hence, $(X,\tau)$ belongs to the class $G_1$.
\end{proof}

\begin{corollary} 
Let $(X,\tau)$ belongs to the class $G_2$. 
Then $\tau$ arise as a dual topology.
\end{corollary}

\begin{question} 
Is it true that all spaces of $G_3$ arise as duals? 
Or conversely, does there exist a space $(X,\tau)$ belonging to the class
$G_3$ such that there is no topology $\sigma$ with $\sigma^d=\tau$?
\end{question}

\section{Examples}

The following examples are due to J. de~Groot, H.~Herrlich, 
G.~E.~Strecker, E.~Wattel (Examples from 4.1 to 4.4, see \cite{GHSW}) and
B. S.~Burdick (Example 4.5, \cite{Bu}) and we repeat them because of
completeness. 
In the most cases, their listed properties follow directly from the fact
that for $T_1$ spaces the dual operator $d$ and the compactness operator
$\rho$ of \cite{GSW} coincide.

\begin{example}
Every compact $T_2$ space belongs to the class $G_1$.
\end{example}

\begin{example}
Every non-compact Hausdorff $k$-space belongs to the class $G_{2_a}$ but 
does not belong to the class $G_1$.
\end{example}

\begin{example} 
Let $W=\omega_1$ be the first uncountable ordinal with its natural order
topology, ${\mathbb N}$ be the discrete space of all natural numbers, 
$Y=W\times{\mathbb N}$ be their product space. 
Let $a\notin Y$, $b\notin Y$, $a\ne b$. 
We set $X=Y\cup\left\{a,b\right\}$ and define the closed subbase $\Phi$
for the topology $\tau$ of the space $(X,\tau)$: 
A subset $F\subseteq X$ belongs to $\Phi$ if and only if fulfills at 
least one of
the following
three conditions:
\begin{itemize}
\item[(i)] 
$F$ is a compact subset of $Y$.
\item[(ii)] 
There are $\alpha\in W$, $n\in {\mathbb N}$ with 
$F = 
\left\{\left(\beta,n\right) : \beta\ge\alpha\right\} \cup 
\left\{a\right\}$
\item[(iii)] 
There are $\alpha\in W$, $\beta\in W$, $n\in{\mathbb N}$ with
$$
F = 
\left\{\left(\gamma,m\right) : \alpha<\gamma\le\beta, n<m\right\} \cup 
\left\{b\right\}
$$
\end{itemize}
Then $(X,\tau)$ belongs to the class $G_{2_b}$ but does not belong to the
class $G_1$.
\end{example}

\begin{example}
Let $X={\mathbb N} \cup \left\{z\right\}$ for some 
$z\in \beta {\mathbb N} \smallsetminus {\mathbb N}$
with the topology $\tau$ induced from the \v Cech-Stone compactification 
$\beta{\mathbb N}$ of ${\mathbb N}$. 
Then $(X,\tau)$ belongs to the class $G_3$ but does not belong to the
class $G_2$.
\end{example}

\begin{example}
Let $(X,\tau)$ be the first uncountable ordinal $X=\omega_1$ with the
topology 
$$
\tau = 
\left\{\left[0,\alpha\right)\smallsetminus F : 
0 \le \alpha \le \omega_1, \text{F is finite}\right\}.
$$ 
Then $(X,\tau)$ belongs to the class $G_4$ but does not belong to the
class $G_3$ (since $\tau^d$ is cocountable, $\tau^{dd}$ is cofinite and
$\tau^{ddd}$ is discrete).
\end{example}

\providecommand{\bysame}{\leavevmode\hbox to3em{\hrulefill}\thinspace}
\providecommand{\MR}{\relax\ifhmode\unskip\space\fi MR }
\providecommand{\MRhref}[2]{%
  \href{http://www.ams.org/mathscinet-getitem?mr=#1}{#2}
}
\providecommand{\href}[2]{#2}

\end{document}